\documentclass[a4paper]{amsart}

\usepackage[dvips]{color}
\usepackage{amssymb}
\usepackage{epsfig}
\usepackage{amsmath}
\usepackage{amsfonts}

\usepackage[english]{babel}
\usepackage{latexsym}

\newcommand{\dd}{\mathrm{d}}

\newcommand*{\DOI}[1]{[DOI:#1]}

\begin{document}

\title[maximal Newtonian resistance]{On the
two-dimensional rotational body\\of maximal
Newtonian resistance\footnote{This is a preprint version 
of the paper published in
\emph{J. Math. Sci. (N. Y.)} {\bf 161}, no.~6, 2009, 811--819.
\DOI{10.1007/s10958-009-9602-0}}}




\author[P. D. F. Gouveia, A. Yu. Plakhov,
D. F. M. Torres]{Paulo D. F. Gouveia,
Alexander Plakhov, Delfim F. M. Torres}

\address{Bragan\c{c}a Polytechnic Institute,
5301-854 Bragan\c{c}a, Portugal}
\email{pgouveia@ipb.pt}

\address{University of Wales, Aberystwyth SY23 3BZ, Ceredigion, UK \newline
\indent (on leave from University of Aveiro, 3810-193 Aveiro, Portugal)}

\email{axp@aber.ac.uk, plakhov@mat.ua.pt}

\address{University of Aveiro, 3810-193 Aveiro, Portugal}
\email{delfim@ua.pt}


\thanks{This work was supported by the
\emph{Centre for Research on Optimization and Control} (CEOC)
from the \emph{Portuguese Foundation for Science and Technology} (FCT),
cofinanced by the \emph{European Community Fund} (ECF)
\textsf{FEDER/POCI 2010}; and by the FCT research project
\textsf{PTDC/MAT/72840/2006}. The first author was also
supported by the ECF \textsf{PRODEP III/5.3/2003}.}


\begin{abstract}
We investigate, by means of computer
simulations, shapes of nonconvex bodies that maximize resistance to their
motion through a rarefied medium, considering that bodies are moving forward
and at the same time slowly rotating. A two-dimensional
geometric shape that confers to the body a resistance
very close to the theoretical supremum value is obtained,
improving previous results.
\end{abstract}

\maketitle


\section{Introduction}

One area of investigation in contemporary mathematics is concerned with the search for shapes of bodies, within predefined classes, which permit the minimization or maximization of the resistance to which they are subjected when they move in rarefied media. The first problem of this nature goes back to the decade of the 1680s, a time when Isaac Newton studied a problem of minimum resistance for a specific class of convex bodies, which moved in media of infinitesimal particles, rarefied to such a degree that it was possible to discount any interaction between the particles, and in which the interaction of these with the bodies could be described as perfectly elastic collisions \cite{newton1686}. More recently we have witnessed important developments in this area with the broadening of study to new classes of bodies and to media with characteristics which are less restrictive. However, almost without exception, the studies which have been published have given special attention to classes of convex bodies ---~the convexity of a body is a sufficient condition for the resistance to be solely a function of singular collisions. Even the various studies of classes of nonconvex bodies which have emerged, especially in the last decade, are based almost always on conditions that guarantee a single impact per particle ---~\cite{brock96,buttazzo93,comte01,robert01}. Only very recently have there begun to emerge some studies supposing multiple reflections, as is the case of the work of Plakhov~\cite{Plakhov03b,plakhov03,Plakhov04}.

In the class of convex bodies, the problem is normally reduced to the minimization of Newton's functional --- an analytical formula for the value of the resistance. But, in the context of nonconvex bodies, there is not any simple formula known for the calculation of the resistance. Even if it is extremely complex, in general, to deal analytically with problems of multiple collisions, for some specific problems of minimization the job has not been revealed to be particularly difficult, there even being some results already available~\cite{Plakhov03b,plakhov03}. If, on the other hand, we consider the problem of maximization, then in this case the solution becomes trivial --- for any dimension, it is enough that the front part of the body is orthogonal to the direction of the movement.

And what if the body exhibits, besides its translational movement, a slight rotational movement? When we think of this kind of problem, we have in mind, for example, artificial satellites, of relatively low orbits, which do not possess any control system which could stabilize their orientation, or other devices in similar conditions. In this situation we imagine that, over its path, the device rotates slowly around itself.

The problem of resistance minimization for rotating nonconvex two-dimensional bodies has already been studied in~\cite{Plakhov04}: it was shown that the reduction of resistance, as compared with the convex
case, does not exceed $1.22\%$.
In its turn, the problem of maximization of the average resistance of bodies in rotation is far from being trivial, in contrast with that which occurs when we deal with purely translational movement.
This class of problems was, therefore, the object of study of the work carried out by the authors in \cite{Plakhov07:CM, Plakhov07}: nonconvex shapes of bodies were investigated which would maximize the resistance that they would have to confront if they moved in rarefied media, and, simultaneously, exhibited a slight rotational movement. With the numerical study which was executed, various geometrical shapes were found which conferred on the bodies rather interesting values of resistance; but it was in later investigations, performed in the follow-up of this work, that the authors managed to arrive at the best of the results --- a two-dimensional shape which confers on the body a resistance very near to its maximum theoretical limit. It is this latest result which now is presented here.


\section{Definition of the problem for the two-dimensional case}
\label{sec:defBidimens}
Consider a disc in slow and uniform rotation, moving in a direction parallel to its plane. We will designate the disc of radius $r$ by $C_r$ and its boundary by $\partial C_r$. We then remove small pieces of the disc along its perimeter, in an $\varepsilon$-neighborhood of $\partial C_r$, with $\varepsilon \in \mathbb{R}_+$ of value arbitrarily small when compared with the value of $r$. We are thus left with a new body $B$ defined by a subset of $C_r$ and characterized by a certain roughness along all its perimeter. The essential question which we put is the following: up to what point can the resistance of a body $B$ be augmented? More than getting to know the absolute value of this resistance, we are principally interested in learning what is the increase which can be obtained in relation to the smooth body (a perfectly circular contour, in this case), that is, learning the normalized value
\begin{equation}
\label{eq:normalizacao}
R(B)=\frac{\text{Resistance}(B)}{\text{Resistance}(C_r)} \text{.}
\end{equation}
It is possible, from the beginning, to know some important reference values for the normalized resistance: $R(C_r)=1$ and the value of the resistance $R(B)$ will have to be found between $0.9878$ (\cite{Plakhov04}) and $1.5$. The value $1.5$ will be hypothetically achieved if all the particles are reflected by the body with the velocity $\mathbf{v}^+$ (velocity with which the particles separate definitively from the body) opposite to the velocity of incidence $\mathbf{v}$ (velocity with which the particles strike the body for the first time), $\mathbf{v}^+=-\mathbf{v}$, a situation in which the maximum momentum is transmitted to the body. It is also possible for us to know the resistance value of some elementary bodies of the type $B$. This is the case, for example, of discs with the contour entirely formed by rectangular indentations which are arbitrarily small or with the shape of rectangular isosceles triangles. As was demonstrated in \cite{Plakhov07}, these bodies are associated with resistances, respectively, of $R=1.25$ and $R=\sqrt{2}$.

Apart from being defined in the disc $C_r$, it is assumed that the body to be maximized is a connected set $B\in\mathbb{R}^2$, with piecewise smooth boundary $\partial B$. Therefore, let us consider a billiard in $\mathbb{R}^2\setminus B$. An infinitesimal particle moves freely, until, upon colliding with the body $B$, it suffers various reflections (one at least) at regular points of its boundary $\partial B$, ending up by resuming free movement which separates it definitively from the body.
Denote by $\text{conv}B$ the convex hull of $B$.
The particle intercepts the $\partial(\text{conv}B)$ contour twice: when it enters into the set $\text{conv}B$ and in the moment that it leaves. $L=|\partial(\text{conv}B)|$ is considered the total length of the curve $\partial(\text{conv}B)$ and the velocity of the particle is in the first and second moments of interception is represented by $\mathbf{v}$ and $\mathbf{v}^+$, and $x$ and $x^+$ the respective points where they occur. As well, the angles which the vectors $-\mathbf{v}$ and $\mathbf{v}^+$ make with the outer normal vector to the section of $\partial(\text{conv}B)$ between the points $x$ and $x^+$ are designated $\varphi$ and $\varphi^+$.
They will be positive if they are defined in the anti-clockwise direction from the normal vector, and negative in the opposite case. With these directions, both $\varphi$ as well as $\varphi^+$ take values in the interval $[-\pi/2,\pi/2]$.

Representing the cavities which characterize the contour of $B$ by subsets $\Omega_1,$ $\Omega_2, \ldots$, which in their total make up the set $\text{conv}B \setminus B$, the normalized resistance of the body $B$ (equation~\eqref{eq:normalizacao}) takes the following form (cf. \cite{Plakhov07}):
\begin{equation}
\label{eq:RB2}
R(B)
=\frac{|\partial(\text{conv}B)|}{|\partial C_r|} \left(\frac{L_0}{L}+\sum_{i \ne 0}{\frac{L_i}{L}R(\tilde{\Omega}_i)}\right)
\text{,}
\end{equation}
being $L_0=|\partial(\text{conv}B)\cap \partial B|$ the length of the convex part of the contour $\partial B$, $L_i=|\partial(\text{conv}B)\cap {\Omega}_i|$, with $i=1,2,\ldots$, the size of the opening of the cavity ${\Omega}_i$, and $R(\tilde{\Omega}_i)$ the resistance of the normalized cavity $\tilde{\Omega}_i$, in relation to a smooth segment of unitary size, with
\begin{equation}
\label{eq:R}
R(\tilde{\Omega}_i)=\frac{3}{8} \int_{-1/2}^{1/2}\int_{-\pi/2}^{\pi/2} \left(
1+\cos\left( \varphi^+(x,\varphi) -\varphi \right)
\right) \cos \varphi\, \dd\varphi\,\dd x
\text{.}
\end{equation}
The function $\varphi^+$ should be seen as the angle of departure of a particle which interacts with a cavity $\tilde{\Omega}_i$ that has opening of unit size and is similar to $\Omega_i$, with the similarity factor $1/L_i$.

From equation~\eqref{eq:RB2}, we understand that the resistance of $B$ can be seen as a weighted mean ($\sum_i L_i/L=1$) of the resistances of the individual cavities which characterize all its boundary (including resistance of the convex part of the boundary), multiplied by a factor which relates the perimeters of the bodies $\text{conv}B$ and $C_r$. Thus, maximizing the resistance of the $B$ body amounts to maximizing the perimeter of $\text{conv}B$ ($|\partial(\text{conv}B)|\le |\partial C_r|$) and the individual resistances of the cavities $\Omega_i$.

Having found the optimal shape $\Omega^*$, which maximizes the functional \eqref{eq:R}, the body of maximum resistance $B$ will be that whose boundary is formed only by the concatenation of small cavities with this shape. We can therefore restrict our problem to the sub-class of bodies $B$ which have their boundary integrally covered with equal cavities, and in doing so admit, without any loss of generality, that each cavity $\Omega_i$ occupies the place of a circle arc of size $\varepsilon\ll r$. As with $L_i=2r\sin(\varepsilon/2r)$, the ratio between the perimeters takes the value
\begin{equation}
\label{eq:RxPerimetros}
\frac{|\partial(\text{conv}B)|}{|\partial C_r|}
=\frac{ \sin(\varepsilon/2r)}{\varepsilon/2r}
\approx 1- \frac{(\varepsilon/r)^2}{24}\text{,}
\end{equation}
or that is, given a body $B$ of a boundary formed by cavities similar to $\Omega$, from~\eqref{eq:RB2} and~\eqref{eq:RxPerimetros}, we conclude that the total resistance of the body will be equal to the resistance of the individual cavity $\Omega$, less a small fraction of this value, which can be neglected when $\varepsilon\ll r$,
\begin{equation}
\label{eq:Raprox}
R(B)\approx R(\Omega)-\frac{(\varepsilon/r)^2}{24}R(\Omega)\text{.}
\end{equation}

Thus, our research has as its objective the finding of cavity shapes $\Omega$ which maximize the value of the functional~\eqref{eq:R}, whose limit we know to be found in the interval
\begin{equation}
\label{eq:RminMax}
1\le \text{sup}_\Omega R(\Omega) \le 1.5
\text{,}
\end{equation}
as is easily proven using~\eqref{eq:R}:
if $\Omega$ is a smooth segment, $\varphi^+(x,\varphi)=-\varphi$ and
$R(\Omega)=\frac{3}{8} \int_{-1/2}^{1/2}\int_{-\pi/2}^{\pi/2} \left(
1+\cos\left(  2\varphi \right)
\right) \cos \varphi\, \dd\varphi\,\dd x=1
$;
in the conditions of maximum resistance $\varphi^+(x,\varphi)=\varphi$,
thus $
R(\Omega)\le\frac{3}{8} \int_{-1/2}^{1/2}\int_{-\pi/2}^{\pi/2} 2 \cos \varphi\, \dd\varphi\,\dd x
=1.5$.


\section{Numerical study of the problem}
\label{sec:estNumer}

In the class of problems which we are studying, only for some shapes of $\Omega$ which are very elementary is it possible to derive an analytical formula of their resistance~\eqref{eq:R}, as we saw in the rectangular and triangular shapes previously referred to. For somewhat more elaborate shapes, the analytical calculation becomes rapidly too complex, if not impossible, given the great difficulty in knowing the function $\varphi^+: [-1/2,1/2]\times[-\pi/2,\pi/2]\rightarrow[-\pi/2,\pi/2]$, which as we know, is intimately related to the format of the cavity $\Omega$. Therefore, recourse to numerical computation emerges as the natural and inevitable approach in order to be able to investigate this class of problems.

There have been developed various computational models which simulate the dynamics of billiard in the cavity. The algorithms of construction of these models, as well as the those responsible for the numerical calculation of the associated resistance, were implemented using the programming language \textsf{C}, given the computational effort involved. The efficiency of the object code, generated by the compilers of \textsf{C}, allowed the numerical approximation of~\eqref{eq:R} to be made with a sufficiently elevated number of subdivisions of the intervals of integration --- between some hundreds and various thousands (up to $5000$). The results were, because of this, obtained with a precision which reached in some cases $10^{-6}$. This precision was controlled by observation of the difference between successive approximations of the resistance $R$ which were obtained with the augmentation of the number of subdivisions.

For the maximization for the resistance of the idealized models, there were used the global algorithms of optimization of the \textit{toolbox} ``\emph{Genetic Algorithm and Direct Search}''  (version 2.0.1 (R2006a), documented in~\cite{toolbox}), a collection of functions which extends the optimization capacities of the MATLAB numerical computation system.
The option for Genetic and Direct search methods is essentially owed to the fact that these do not require any information about the gradient of the objective function nor about derivatives of a higher order --- as the analytical form of the resistance function is in general unknown (given that it depends on $\varphi^+(x,\varphi)$), this type of information, if it were necessary, would have to be obtained by numerical approximation, something which would greatly impede the optimization process. The MATLAB computation system (version 7.2 (R2006a)) was also chosen because it had functionalities which allowed it to be used for the objective function the subroutine compiled in \textsf{C} of resistance calculation, as well as the $\varphi^+(x,\varphi)$ function invoked in itself.


\subsection{``Double Parabola'': a two-dimensional shape which maximizes resistance}
\label{sec:DuplaParab}

In the numerical study which the authors carried out in \cite{Plakhov07:CM, Plakhov07}, shapes of $\Omega_f$ defined by continuous and piecewise differentiable $f:[-1/2,1/2]\rightarrow \mathbb{R}_+$ functions were sought for:
\begin{equation}
\label{eq:Omegaf}
\Omega_f=\left\{(x,y):\,-1/2\le x \le 1/2,\; 0\le y \le f(x)\right\}\text{,}
\end{equation}
with the interval $[-1/2, 1/2] \times \{ 0 \}$ being the opening.

The search for the maximum resistance was begun in the class of continuous functions $f$ with derivative $f'$ piecewise constant, broadening later to the study of classes of functions with the second derivative $f''$ piecewise constant. Not having been able with these shapes to exceed the value of resistance $R=1.44772$, we decided, in this new study, to extend the search to shapes different from those considered in~\eqref{eq:Omegaf}. We studied shapes $\Omega^g$ defined by functions $x$ of $y$ of the following form:
\begin{equation}
\label{eq:Omegag}
\Omega^g=\left\{(x,y):\,0\le y \le h,\; -g(y)\le x \le g(y)\right\}\text{,}
\end{equation}
where $h>0$ and $g:[0,h] \rightarrow \mathbb{R}_0^+$ is a continuous function with $g(0)=1/2$ and $g(h)=0$.

Similarly to the study which was carried out for the sets $\Omega_f$, in the search for shapes $\Omega^g$, the functions $g$ were considered piecewise linear and piecewise quadratic. If in the classes of linear functions it was not possible to achieve a gain in resistance relative to the results obtained for the sets $\Omega_f$, in the quadratic functions the results exceeded the highest expectations: there was found a shape of cavity $\Omega^g$ which presented the resistance $R=1.4965$, a value very close to its theoretical limit of $1.5$.
There were also carried out some tests with polynomial functions of higher order or described by specific conical sections, but, not having verified any additional gain in the maximization of resistance, it was decided not to report the respective results. There therefore follows the description of the best result which was obtained, encountered in the class of quadratic functions $x=\pm g(y)$.

\begin{figure}[!ht]
\begin{center}
\begin{tabular}{c c}
\includegraphics*[height=0.36\columnwidth]{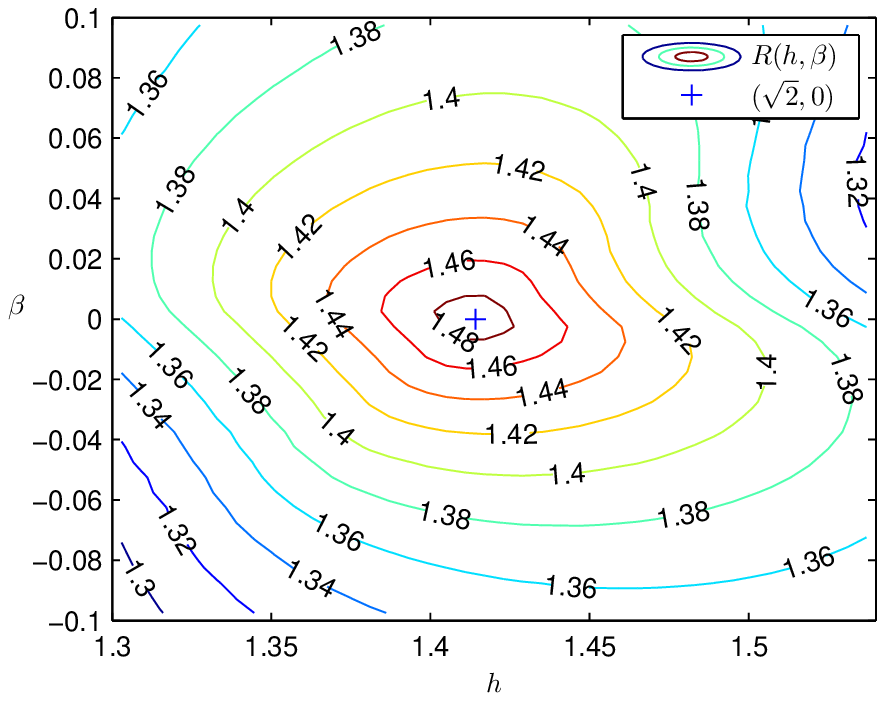}
&
\includegraphics*[height=0.36\columnwidth]{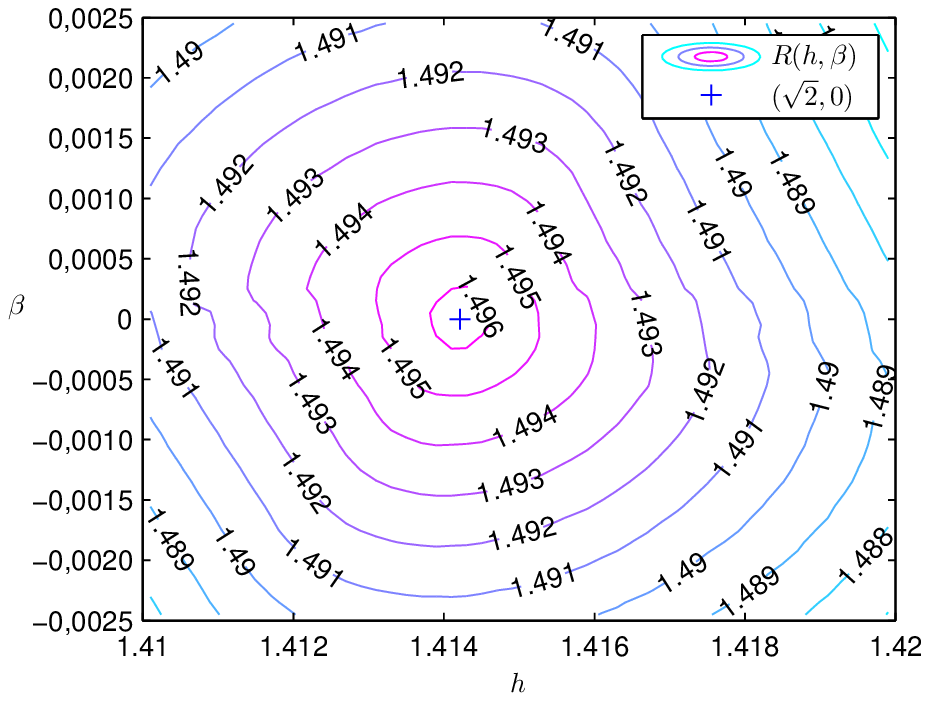} \\
(a)&(b)\\
\end{tabular}
\caption{Level curves of the $R(h,\beta)$ function.}
\label{fig:curvNivel}
\end{center}
\end{figure}

The value of resistance of the sets $\Omega^g$ were studied, just as defined in~\eqref{eq:Omegag}, in the class of quadratic functions
$$
g_{h,\beta}(y) =
 \alpha y^2 + \beta y +1/2, \text{ para } 0 \le y \le h\,,
$$
where $h>0$ and $\alpha= \frac{-\beta h -1/2}{h^2}$ (given that $g_{h,\beta}(h) =0$). In the optimization of the curve, the two parameters of the configuration were made to vary: $h$, the height of the $\partial\Omega^g$ curve, and $\beta$, in its slope in the origin ($g'(0)$). In this class of functions the algorithms of optimization converge rapidly towards a very interesting result:
the maximum resistance was reached with $h=1.4142$ and $\beta=0.0000$, and assumed the value $R=1.4965$, that is, a value $49.65\%$ above the resistance of the rectilinear segment.
This result seems to us really interesting:
\begin{enumerate}
\item it represents a considerable gain in the value of the resistance, relative to the best result obtained earlier (in \cite{Plakhov07:CM, Plakhov07}),
which was situated $44.77\%$ above the reference value;
\item the corresponding set $\Omega^g$ has a much more simple shape than that of set $\Omega_f$ associated with the best earlier result, since it is formed by two arcs of symmetrical parabolas, while the earlier one was made up of fourteen of these arcs;
\item this new resistance value is very near to its maximal theoretical limit, which, as is known, is found $50\%$ above the value of reference;
\item the optimal parameters appear to assume value which give to the set $\Omega^g$ a configuration with very special characteristics, as in what follows will be understood.
\end{enumerate}
Note that the optimal parameters appear to approximate the exact values $h=\sqrt{2}$ and $\beta=0$.
The graphical representation of the function $R(h,\beta)$ through the level curves, figure~\ref{fig:curvNivel}, are effectively in concordance with this possibility
 --- note that the level curves appear perfectly centered on the $(\sqrt{2},0)$ coordinates; marked on the figure by ``$+$'' .

\begin{figure}[!ht]
\begin{center}
\includegraphics*[width=0.45\columnwidth]{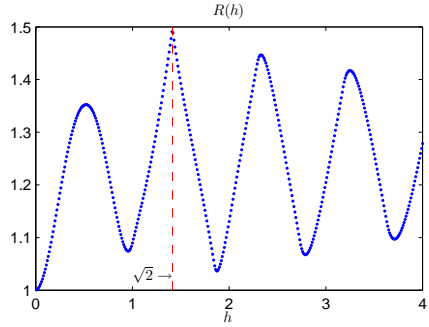}
\caption{Resistance graphic $R(h)$ for $\beta=0$.}
\label{fig:Res_h}
\end{center}
\end{figure}

Note also, in figure~\ref{fig:Res_h}, the resistance graph $R(h)$ for $\beta=0$, where it can equally be perceived that there is a surprising elevation of resistance when $h\rightarrow\sqrt{2}$. Thus the resistance of the $\Omega^{g_{h\beta}}$ cavity was numerically calculated with the exact values $h=\sqrt{2}$ and $\beta=0$, the result having confirmed the value $1.49650$.

\begin{figure}[!ht]
\begin{center}
\begin{tabular}{c c c}
\includegraphics*[height=0.3\columnwidth]{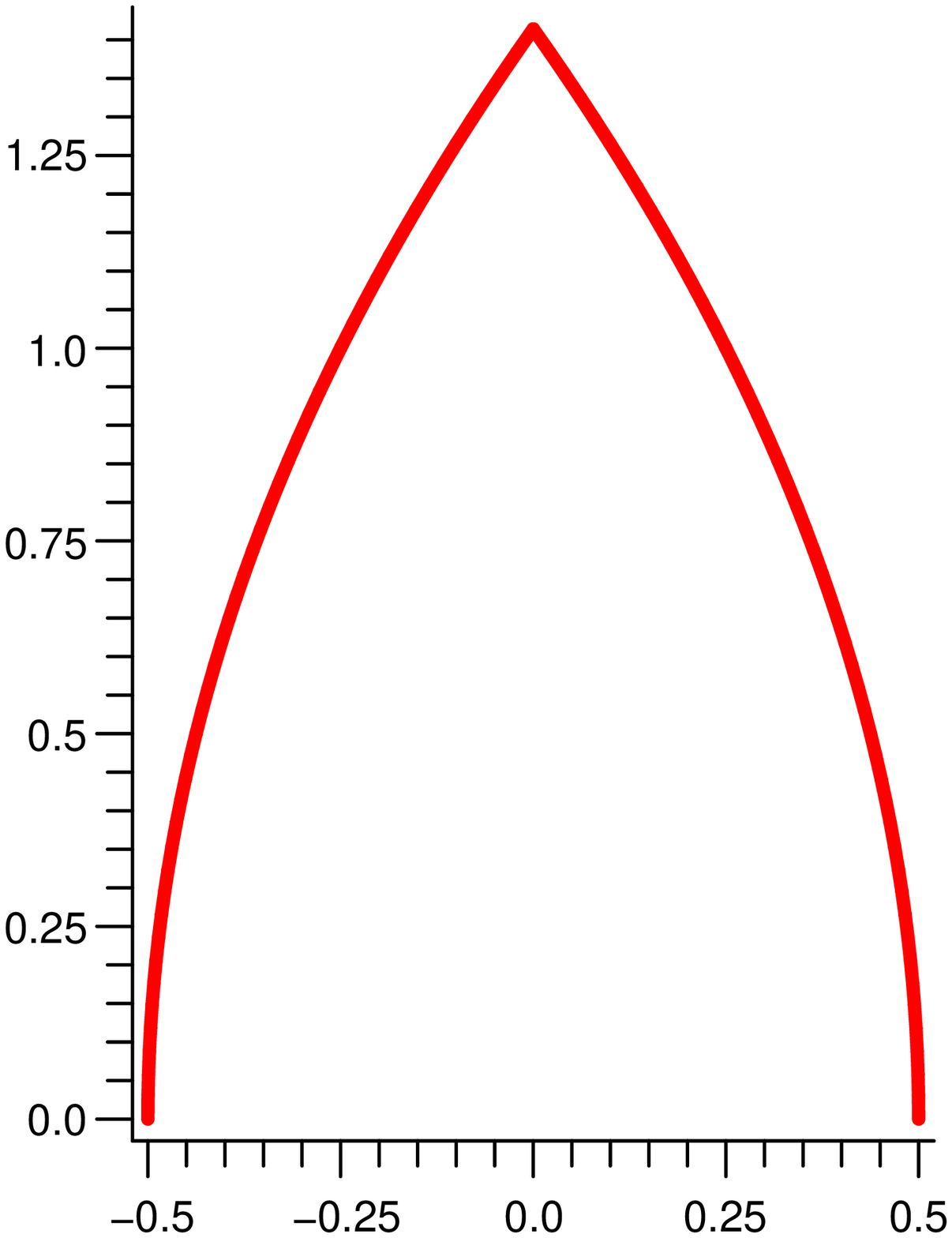}
&
&
\includegraphics*[height=0.3\columnwidth]{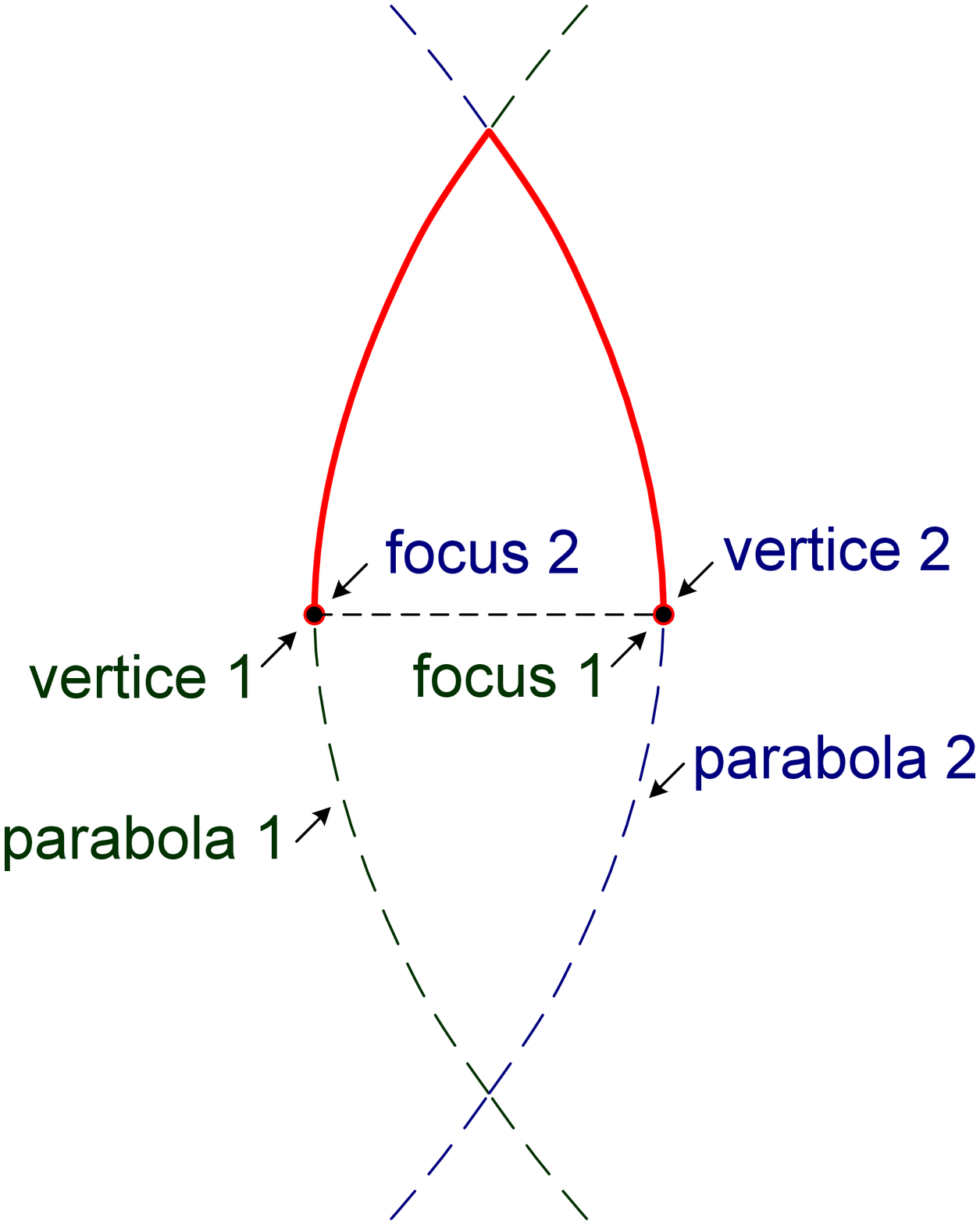}
 \\
(a)&&(b)\\
\end{tabular}
\caption{(quasi) Optimal 2D shape --- the \emph{Double Parabola}.}
\label{fig:parabOpt}
\end{center}
\end{figure}

The shape of the set $\Omega^{g_{h,\beta}}$ with $h=\sqrt{2}$ and $\beta=0$ is a particular case with which is associated special characteristics which could justify the elevated value of resistance presented. The two sections of the shape are similar arcs of two parabolas with the common horizontal axis and concavities turned one towards the other --- see figure~\ref{fig:parabOpt}. But the particularity of the configuration resides in the fact that the axis of the parabolas coincides with the line of entry of the cavity (axis of $x$), and that the focus of each one coincides with the vertex of the other.

In order to be easily referred to, this shape of cavity (figure~\ref{fig:parabOpt}a) will be, from here on, named simply ``\emph{Double Parabola}''. Thus, in the context of this paper, the term ``Double Parabola''  should be always understood as the name of the cavity whose shape is described by two parabolas which, apart from being geometrically equal, find themselves ``nested'' in the particular position to which we have referred.

Since the resistance of the Double Parabolas assumes a value which is very close to its theoretical limit, in a final attempt to achieve this limit, it was resolved to extend the study even further to other classes of functions $g(y)$ which admit the Double Parabola as a particular case or which allow proximate configurations of this nearly optimal shape. In all these cases the best results were invariably obtained when the shape of the curves approximated the shape of the Double Parabola, without ever having overtaken the value $R=1.4965$.


\section{Characterization of the reflections in the ``Double Parabola''}
\label{sec:caract}

\begin{figure}[!hb]
\begin{center}
\begin{tabular}{c c c}
\includegraphics*[width=0.25\columnwidth]{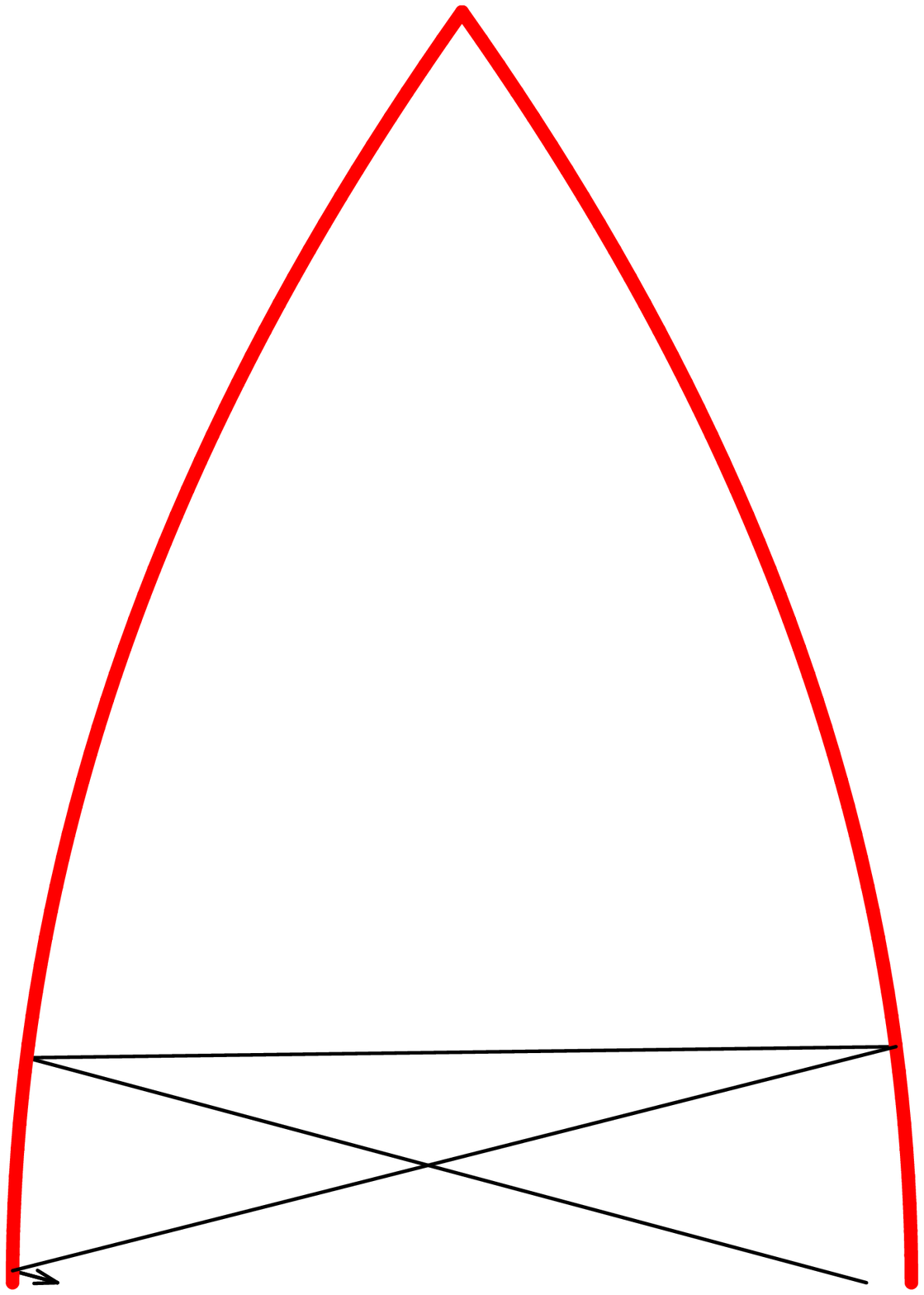} &
\includegraphics*[width=0.25\columnwidth]{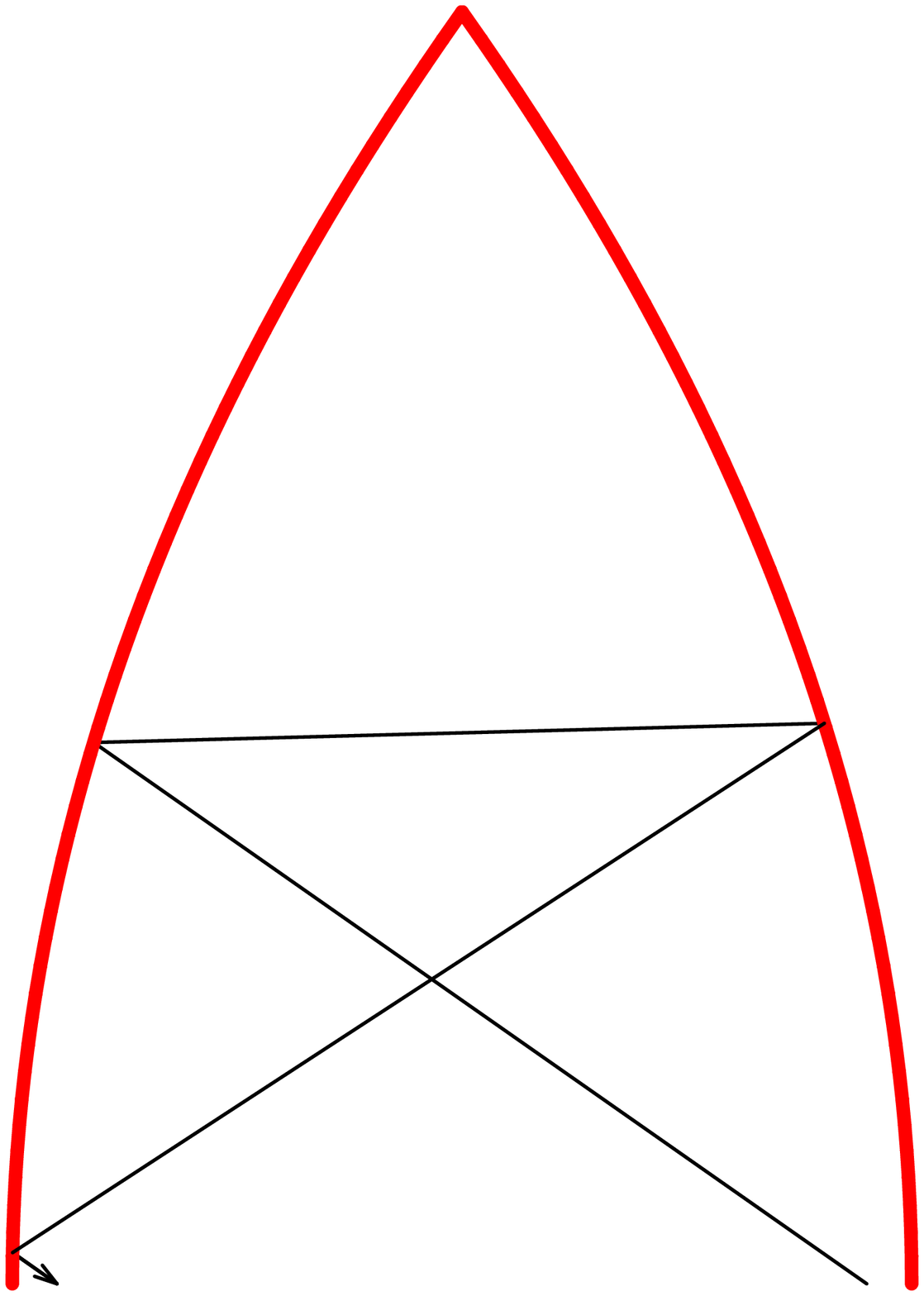} &
\includegraphics*[width=0.25\columnwidth]{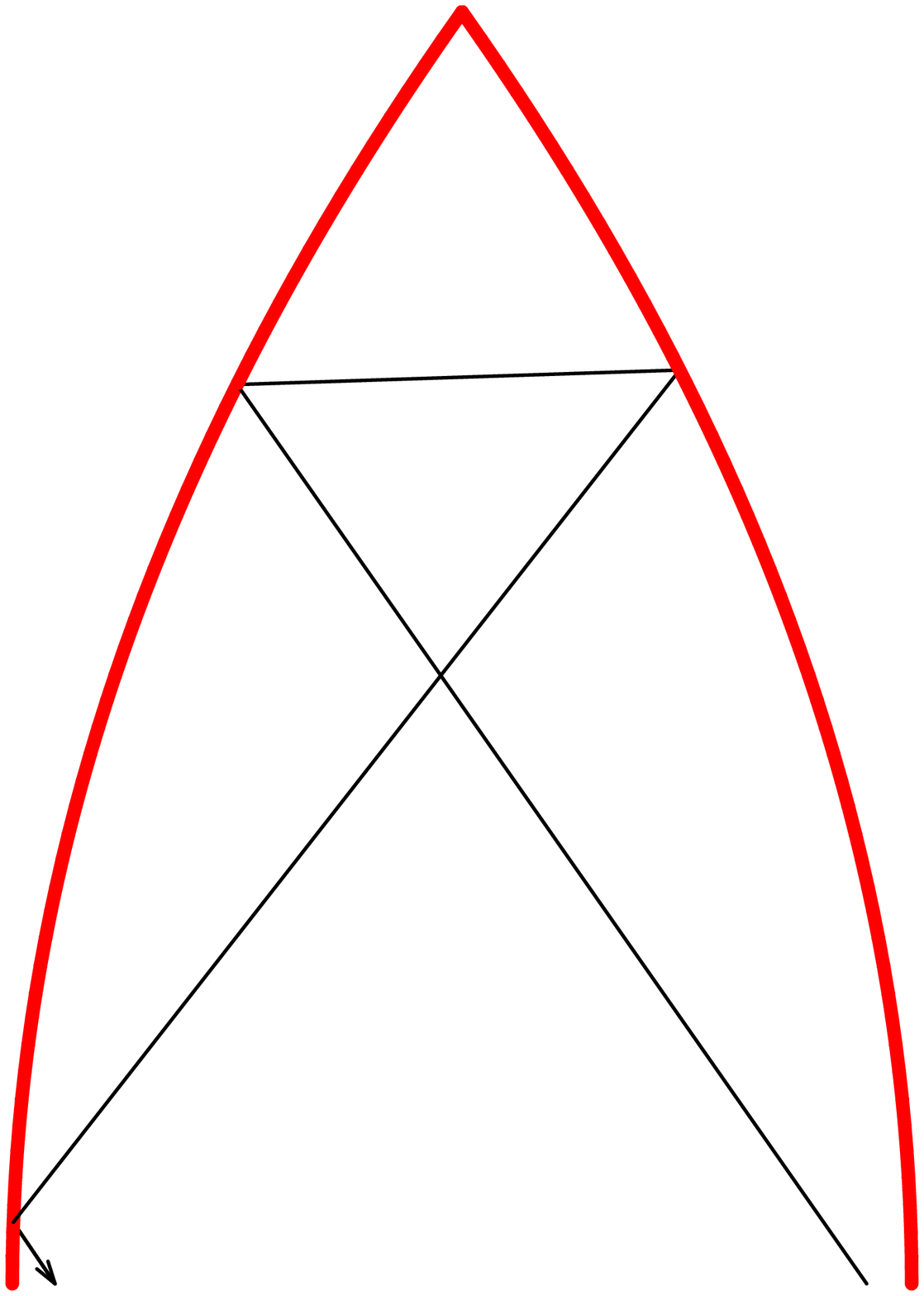} \\
(a) $x=0.45$, $\varphi=75^\circ$.&(b) $x=0.45$, $\varphi=55^\circ$.&(c) $x=0.45$, $\varphi=35^\circ$.\\
\includegraphics*[width=0.25\columnwidth]{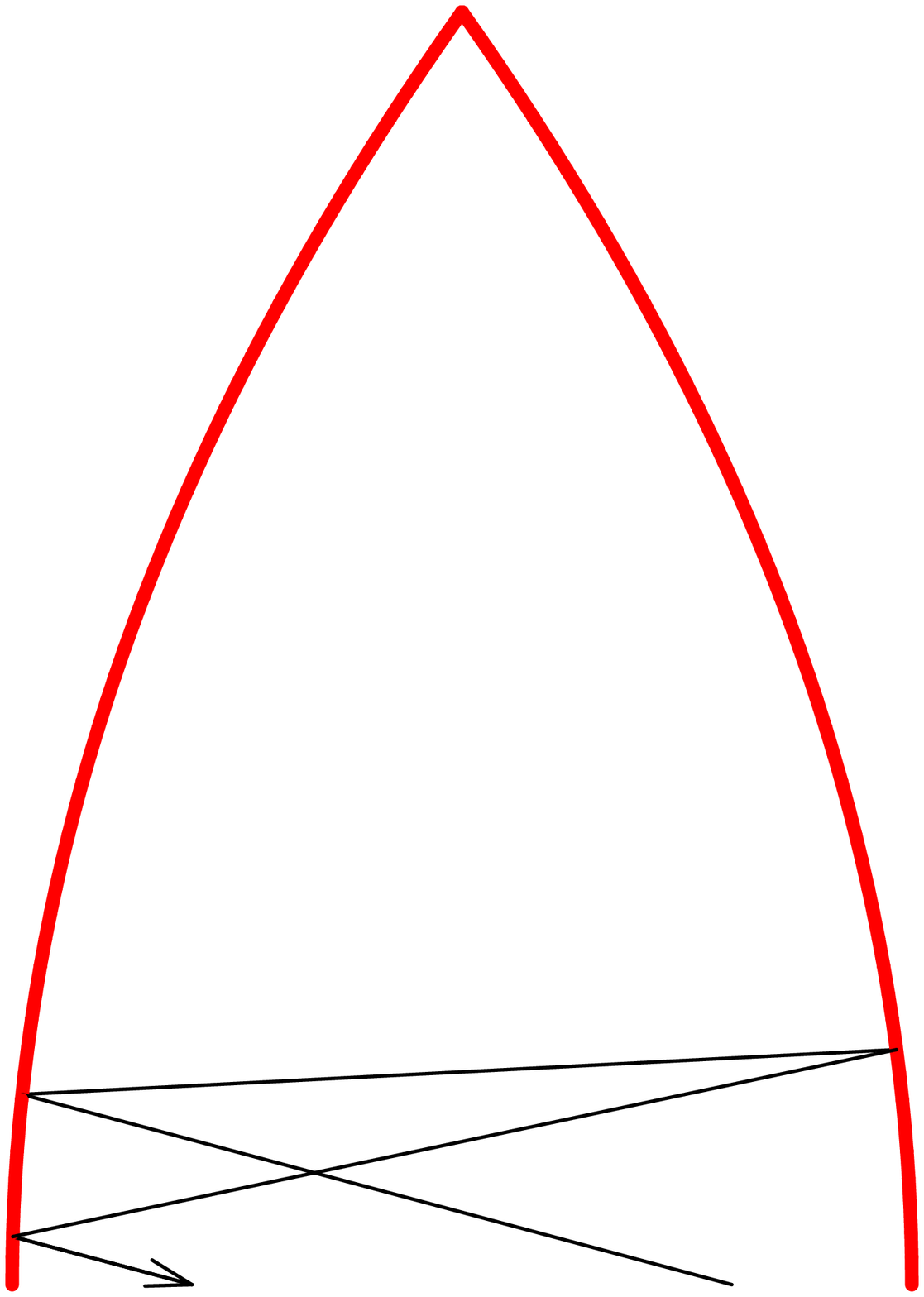} &
\includegraphics*[width=0.25\columnwidth]{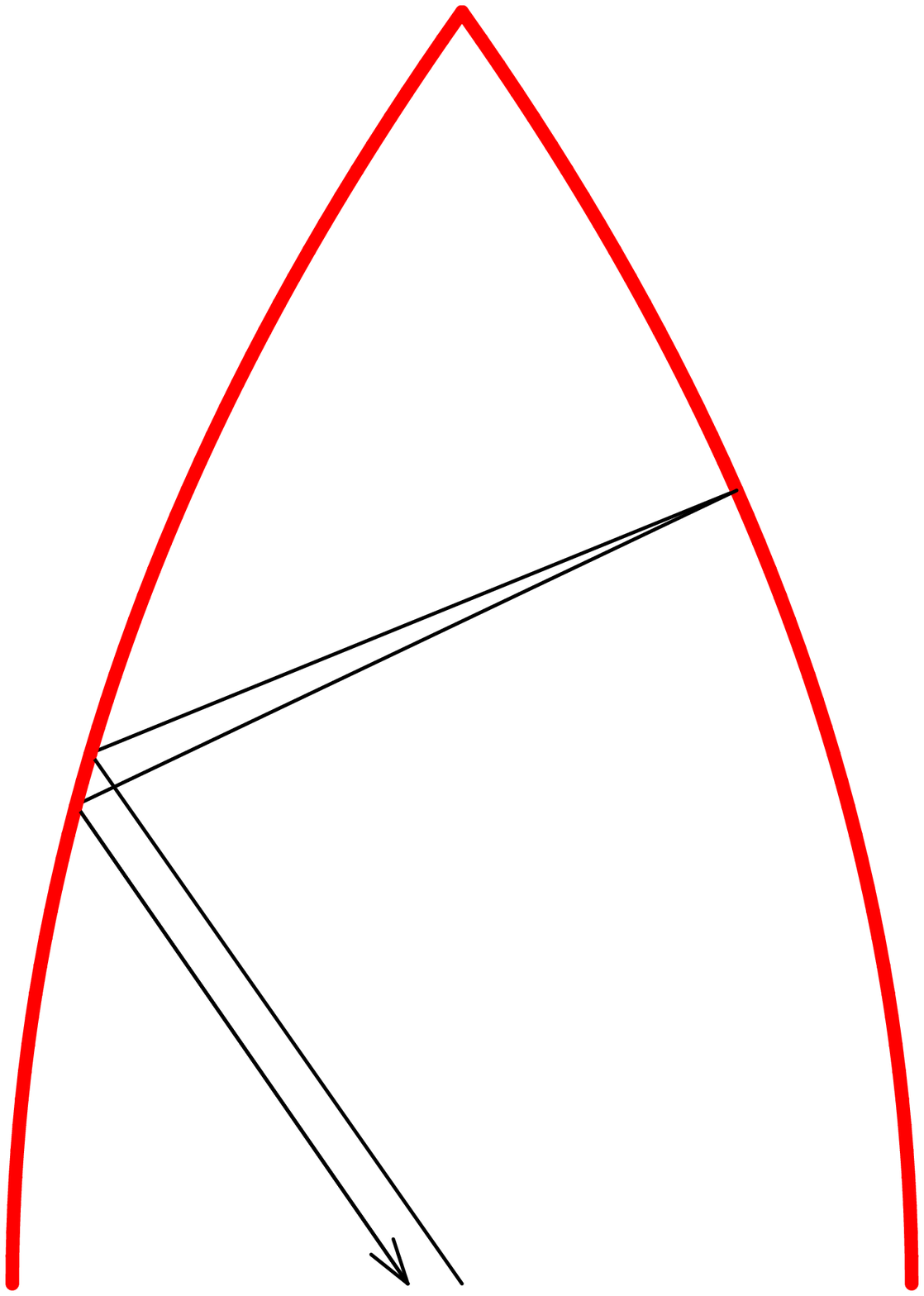} &
\includegraphics*[width=0.25\columnwidth]{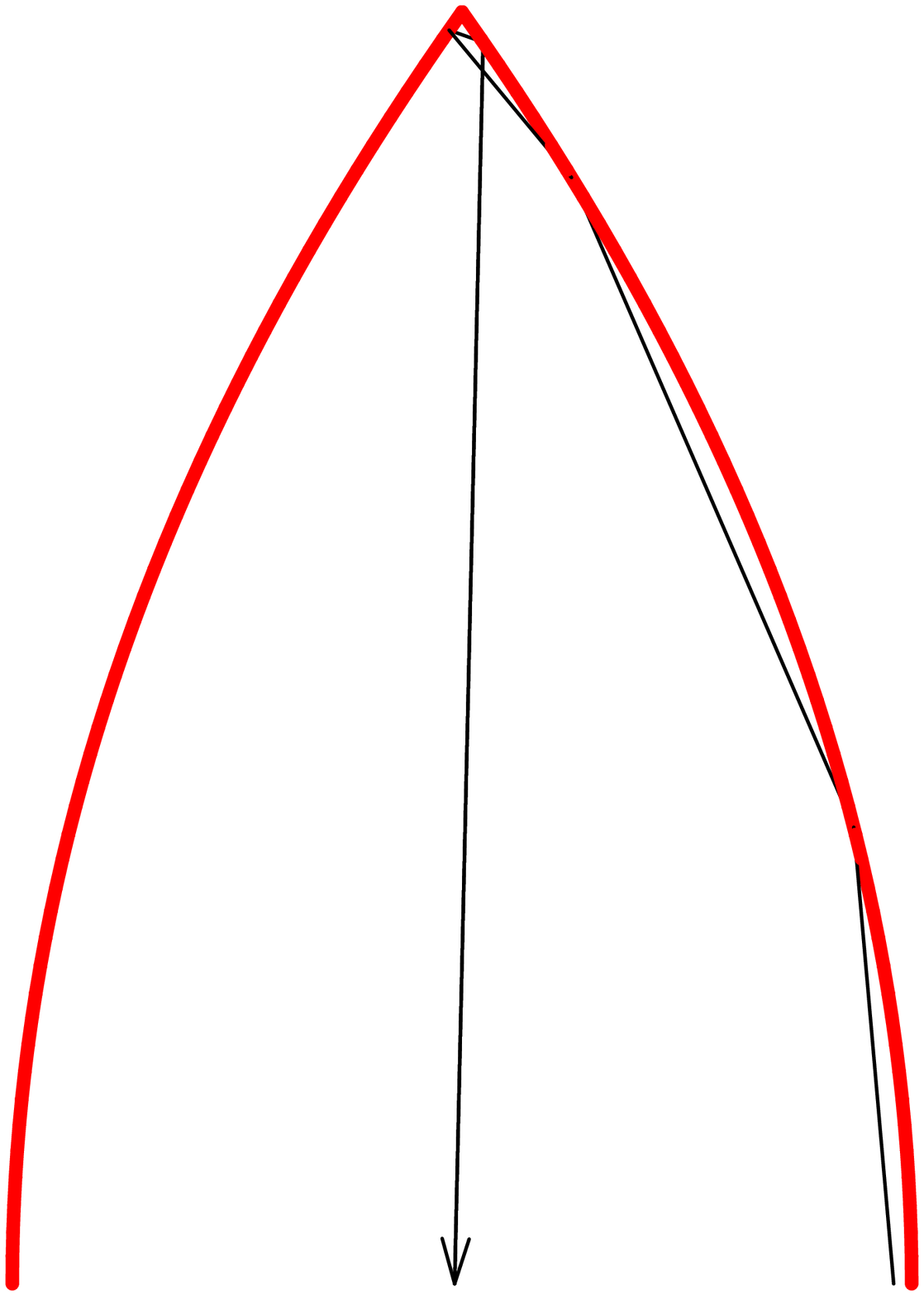} \\
(d) $x=0.3$, $\varphi=75^\circ$.&(e) $x=0.0$, $\varphi=35^\circ$.&(f) $x=0.48$, $\varphi=5^\circ$. \\
\end{tabular}
\caption{Example of trajectories obtained with the computational model.}
\label{fig:trajectorias}
\end{center}
\end{figure}

Each one of the illustrations of figure~\ref{fig:trajectorias} shows, for the ``Double Parabola'', a concrete trajectory, obtained with our computational model.

It is comforting to verify that, with the exception of one trajectory, in all the others the particle emerges from the cavity with a velocity which is nearly opposite to that which was its entry velocity. This is the ``symptom''  which unequivocally characterizes a cavity of optimal performance. Even in the case of the trajectory of the illustration (f), the direction of the exit velocity appears not to vary greatly from that of entry.

If we analyze the five first illustrations, we may verify that there exists something in common in the behavior of the particle: in describing the trajectory, the particle is always subject to three reflections. This appears to be a determinant characteristic for the approximation of the two angles of entry and exit.
Although this conviction is by nature essentially empirical, the results of the analytical study which the first author carried out in \cite{gouvPhD} are heading in the direction of confirming that one very significant part of the ``benign''  trajectories --- those in which the vectors velocity of entry and of exit are nearly parallel; we call them so because they represent positive contributions to the maximization of resistance --- suppose exactly three reflections.
In that study the author managed to demonstrate some important properties which help in the consolidating the numeric results presented for ``Double Parabola''.
In particular, it was shown that
\begin{itemize}
\item There are no trajectories of fewer than $3$ reflections;
\item For angles of entry $\varphi$ outside the interval $(-\varphi_0,\varphi_0)$, with $\varphi_0=\arctan\left(\frac{\sqrt{2}}{4}\right)$ $\simeq 19.47^\circ$, all the trajectories are of $3$ reflections;
\item In trajectories with $4$ or more reflections, the angular difference is delimited by $2\varphi_0$: $\left|\varphi-\varphi^+\right|<2\varphi_0$.
\end{itemize}


\section{Other possible applications}
\label{sec:outrasAplic}

Besides maximizing Newtonian resistance, it is exciting to verify that the potentialities of the Double Parabola shape found by us could also reveal themselves to be very interesting in other areas of practical interest.

Given the characteristics of reflection which the Double Parabola shape presents, we can rapidly conceive for it a natural propensity for being able to be used with success in the design of retroreflectors. Retroreflectors are devices which send light or other incident radiation back to the emitting source. Ideally the retroreflector should undertake this function independently of the angle of incidence, something which does not happen with existing devices.

As we have seen, the Double Parabola, although it does not guarantee the perfect inversion of all the incident radiation, carries out this function with great success:
it guarantees a good approximation of the directions of the incident and reflected flows for a significant part of the angles of incidence, and even for the rest it does not permit that the lag reaches elevated values. We predict, therefore, as very promising its possible utilization in the definition of new geometries for the optical elements which make up retroreflecting surfaces.

Retroreflecting devices, although they may be used in a wide range of technological areas, as is the case for example of optical communications in open space, are used in a massive way in the automobile industry and in roadway signalling, and thus we can understand very easily their utility. In \cite{gouvPhD} is presented a exploratory study about a possible way to take advantage of our result.


\section{Conclusion}
\label{sec:concl}

\begin{figure}[!ht]
\begin{center}
\includegraphics[width=0.30\columnwidth]{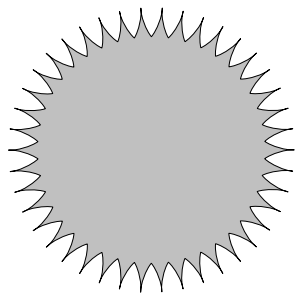}
\caption{(quasi) Optimal 2D body.}
\label{fig:discoOptimo}
\end{center}
\end{figure}

In the continuation of the study carried out previously by the authors in \cite{Plakhov07:CM, Plakhov07}, with the work now presented it has been possible to obtain an original result which appears to us to have great scope: the algorithms of optimization converge for a geometrical shape very close to the ideal shape --- the \emph{Double Parabola}. This concerns a form of roughness which confers a nearly maximal resistance (very close to the theoretical majorant) to a disc which, not only travels in a translational movement but also rotates slowly around itself. In figure~\ref{fig:discoOptimo} one of these bodies is shown.

Noting that the contour of the presented body is integrally formed by $42$ cavities $\Omega$ with the shape of a Double Parabola, each one of which with a relative resistance of $1.49650$, from \eqref{eq:RB2} and \eqref{eq:RxPerimetros} we conclude that $R(B)=\frac{\sin(\pi/42)}{\pi/42}R(\Omega)\approx 1.4951$ is the total resistance of the body, a value $49.51\%$ above the value of resistance of the corresponding disc of smooth contour (the smallest disc which includes the body). We know that if the body were formed by a sufficiently elevated number of these cavities, its resistance would even reach the value $1.4965$, but the example presented is sufficient in order for us to understand how close we are to the $50\%$, the known theoretical majorant.

Taking into account that the majorant $1.5$ is a theoretical result which only
signifies the non-existence of 2D bodies (with boundary formed by identical cavities) that exceed this resistance value, the shape found by us could even embody an optimal solution.
In confirming this hypothesis --- that the value of resistance which characterizes the Double Parabola effectively embodies the maximum limit which it is possible to reach with real two-dimensional shapes --- our result would earn redoubled importance. Although easy developments are not foreseeable, this is an important problem which remains open, waiting for future contributions which, if they do not outdo our result, will permit the reinforcement of our conjecture.



\end{document}